\documentstyle[11pt]{article}
 
\newcommand{\ncm}{\newcommand}
\ncm{\rncm}{\renewcommand}
\rncm{\sec}{\setc{0}\section}
\ncm{\beq}{\begin{equation}}
\ncm{\eeq}{\end{equation}}
\ncm{\bea}{\begin{eqnarray}}
\ncm{\beanon}{\begin{eqnarray*}}
\ncm{\eea}{\end{eqnarray}}
\ncm{\eeanon}{\end{eqnarray*}}
 
\rncm{\theequation}{\thesection.\arabic{equation}}
\ncm{\setc}[1]{\setcounter{equation}{#1}}
\newcounter{eqnr}
\newenvironment{eqnarrayabc}{\stepcounter{equation}
  \setcounter{eqnr}{\value{equation}}\setc{0}
\rncm{\theequation}{\thesection.\arabic{eqnr}\alph{equation}}
  \begin{eqnarray}}{\end{eqnarray}\setc{\value{eqnr}}}
 
\ncm{\beabc}{\begin{eqnarrayabc}}
\ncm{\eeabc}{\end{eqnarrayabc}}
\ncm{\nn}{\nonumber \\}
 
\newtheorem{thm}{Theorem}[section]
\newtheorem{prop}[thm]{Proposition}
\newtheorem{lem}[thm]{Lemma}

\newtheorem{defi}[thm]{Definition}

\def\id{\hbox{id}\,}

\def\Center{\hbox{Center}\,}
\def\Hil{{\cal H}}
\def\one{{\bf 1}}
\def\two{{\bf 2}}
\def\three{{\bf 3}}
\def\1{{\thinmuskip=5.5mu 1\!1\thinmuskip=3mu}}
\def\du1{{\hat\1}}
\def\o{\otimes}

\def\p{^{\prime}}
\def\inv{^{-1}}
\def\cop{\Delta}
\def\ducop{{\hat\cop}}
\def\eps{\varepsilon}
\def\dueps{{\hat\eps}}
\def\duS{{\hat S}}
\def\duA{{\hat A}}

\def\c{_{(1)}}
\def\cc{_{(2)}}
\def\ccc{_{(3)}}

\def\PL{\sqcap^L}
\def\PR{\sqcap^R}

\def\duPL{{\hat\PL}}
\def\duPR{{\hat\PR}}
\def\cros{\raise1.9pt\hbox{$\scriptscriptstyle
          > $}\!\raise1.5pt\hbox{$\scriptstyle\triangleleft\,$}}

\def\C{\,{\raise1.5pt\hbox{$\scriptscriptstyle |$}
        \thinmuskip=4mu \!\!C\thinmuskip=3mu}}

\def\bra{\langle}
\def\ket{\rangle}
\def\la{\rightharpoonup}
\def\ra{\leftharpoonup}
\def\qed{\hfill {\it Q.e.d.}}
\def\Proof{{\em Proof}\,:\ }
\ncm\amalgo[1]{{\lower9pt\hbox{$\o$}\atop\raise2pt\hbox{
            $\scriptscriptstyle #1$}}}
 
\def\amo{{\lower9pt\hbox{$\otimes$}\atop\raise2pt\hbox{
          $\scriptscriptstyle A^R\equiv\hat A^L$}}}
\def\LH{{\cal L}(\Hil)}
\ncm\eq[1]{\stackrel{(\ref{#1})}{=}}
 
\def\oR{\Box\kern -8pt\raise
   2pt\hbox{$\scriptscriptstyle\swarrow$}\,}
\def\oL{\Box\kern -8.5pt\raise
   3pt\hbox{$\scriptscriptstyle\nwarrow$}\,}
\def\HK{\Hil\amalgo{\psi}{\cal K}\, }
\def\KH{{\cal K} \amalgo {\psi^o} \Hil \,}
\def\half{{1\over 2}}
\def\target{\Hil_{\alpha_{02}}\amalgo{\eps}\ _{\alpha_{01}}\Hil}
\def\source{\Hil_{\alpha_{12}}\amalgo{\eps^o}\ _{\alpha_{02}}\Hil}
\def\trimodule{\parbox[c]{1in}{
  \begin{picture}(30,30)(0,20)
  \put(5,45){\line(1,0){40}}
  \put(5,45){\line(1,-1){20}}
  \put(45,45){\line(-1,-1){20}}
  \put(0,46){$\scriptscriptstyle 0$}
  \put(47,46){$\scriptscriptstyle 2$}
  \put(24,19){$\scriptscriptstyle 1$}
  \end{picture}
}}
\def\sourcepic{\parbox[c]{1in}{
\begin{picture}(50,50)(-5,-10)
\put(7,18){$\scriptstyle \one$}
\put(18,7){$\scriptstyle \two$}
\put(0,0){\line(1,0){30}}
\put(0,0){\line(0,1){30}}
\put(0,0){\line(1,1){30}}
\put(30,0){\line(0,1){30}}
\put(0,30){\line(1,0){30}}
\put(-5,33){$\scriptstyle 0$}
\put(-5,-5){$\scriptstyle 1$}
\put(32,-5){$\scriptstyle 2$}
\put(32,33){$\scriptstyle 3$}
\end{picture}
}}
\def\targetpic{\parbox[c]{1in}{
\begin{picture}(50,50)(-5,-10)
\put(8,7){$\scriptstyle \one$}
\put(18,18){$\scriptstyle \two$}
\put(0,0){\line(1,0){30}}
\put(0,0){\line(0,1){30}}
\put(0,30){\line(1,-1){30}}
\put(30,0){\line(0,1){30}}
\put(0,30){\line(1,0){30}}
\put(-5,33){$\scriptstyle 0$}
\put(-5,-5){$\scriptstyle 1$}
\put(32,-5){$\scriptstyle 2$}
\put(32,33){$\scriptstyle 3$}
\end{picture}
}}
\def\tetra{\parbox[c]{1in}{
\begin{picture}(50,50)(-5,-10)
\put(0,0){\line(1,0){30}}
\put(0,0){\line(0,1){30}}
\put(0,0){\line(1,1){13}}
\put(17,17){\line(1,1){13}}
\put(0,30){\line(1,-1){30}}
\put(30,0){\line(0,1){30}}
\put(0,30){\line(1,0){30}}
\put(-5,33){$\scriptstyle 0$}
\put(-5,-5){$\scriptstyle 1$}
\put(32,-5){$\scriptstyle 2$}
\put(32,33){$\scriptstyle 3$}
\end{picture}
}}
\rncm{\t}[3]{(^{#1#3}_{\,#2})}
\def\penta{\parbox[c]{1in}{
  \setlength{\unitlength}{2pt}
  \begin{picture}(35,35)
  \put(5,15){\line(1,3){5}}       
  \put(5,15){\line(3,-2){15}}     
  \put(5,15){\line(5,3){25}}      
  \put(10,30){\line(1,0){20}}     
  \put(20,5){\line(2,5){10}}      
  \put(20,5){\line(3,2){15}}      
  \put(30,30){\line(1,-3){5}}     
  \put(7,31){$\scriptstyle 0$}
  \put(2,13){$\scriptstyle 1$}
  \put(19,1){$\scriptstyle 2$}
  \put(37,13){$\scriptstyle 3$}
  \put(31,31){$\scriptstyle 4$}
  \put(11,23){$\scriptstyle\one$}
  \put(15,14){$\scriptstyle\two$}
  \put(27,14){$\scriptstyle\three$}
  \end{picture}
}}
 
\def\pentta{\parbox[c]{1in}{
  \setlength{\unitlength}{2pt}
  \begin{picture}(35,35)
  \put(5,15){\line(1,3){5}}       
  \put(5,15){\line(3,-2){15}}     
  \put(5,15){\line(5,3){25}}      
  \put(10,30){\line(1,0){20}}     
  \put(5,15){\line(1,0){30}}      
  \put(20,5){\line(3,2){15}}      
  \put(30,30){\line(1,-3){5}}     
  \put(7,31){$\scriptstyle 0$}
  \put(2,13){$\scriptstyle 1$}
  \put(19,1){$\scriptstyle 2$}
  \put(37,13){$\scriptstyle 3$}
  \put(31,31){$\scriptstyle 4$}
  \put(11,23){$\scriptstyle\one$}
  \put(18,9){$\scriptstyle\two$}
  \put(23,18){$\scriptstyle\three$}
  \end{picture}
}}

\begin{document}
\large
\title{\bf Weak C$^*$-Hopf Algebras and Multiplicative
Isometries\\}
 
\author{\sc Gabriella B\"ohm $^1$,  Korn\'el Szlach\'anyi $^2$\\
\\
Research Institute for Particle and Nuclear Physics, Budapest\\
H-1525 Budapest 114, P.O.B. 49, Hungary}

\date{}
 
\maketitle
 
\footnotetext[1]{E-mail:
  BGABR@rmki.kfki.hu\\
  Supported by the Hungarian Scientific Research Fund, OTKA --
  T 016 233}
 
\footnotetext[2]{E-mail:
  SZLACH@rmki.kfki.hu\\
  Supported by the Hungarian Scientific Research Fund,
  OTKA -- T 020 285.}

\vskip 2truecm

\begin{abstract}
We show how the data of a finite dimensional weak $C^*$-Hopf
algebra can be encoded into a pair $(\Hil,V)$ where $\Hil$ is a
finite dimensional Hilbert space and
$V\colon\Hil\o\Hil\to\Hil\o\Hil$ is a partial isometry satisfying,
among others, the pentagon equation. In case of $V$ being unitary
we recover the Baaj-Skandalis multiplicative unitary of the
discrete compact type. Relation to the pseudomultiplicative unitary
approach proposed by J.-M. Vallin and M. Enock is also discussed.
\end{abstract}

\newpage
 
\section{Introduction}
The fundamental operator in Kac algebra theory \cite{ESch}
or the multiplicative unitary in $C^*$-Hopf algebras \cite{BaSk}
is a unitary operator
$V\colon\Hil\o\Hil\to\Hil\o\Hil$ satisfying the pentagon equation
$V_{23}V_{12}=V_{12}V_{13}V_{23}$ on the three-fold tensor
product of the Hilbert space $\Hil$.
It encodes information about the structure of a quantum group $A$
and its dual $\duA$ in a symmetric way.
If $\Hil$ is finite dimensional then a multiplicative unitary is
the complete information necessary to
determine a unique finite dimensional $C^*$-Hopf algebra
\cite{BaSk}. In the
infinite dimensional case additional assumptions are necessary:
These are the regularity and irreducibility assumptions in the
work of Baaj and Skandalis.
 
If $A$ is a finite dimensional $C^*$-Hopf algebra then a
multiplicative unitary on the Hilbert space of the left regular
representation can be given by the formula $V(x\o y)=x\c\o x\cc y$
where $x\mapsto \cop(x)\equiv x\c\o x\cc$ denotes the coproduct
on $A$ and $x,y\in\Hil\equiv A$. As it has been noticed in
\cite{BSz}
if $A$ is only a weak $C^*$-Hopf algebra then the $V$ defined by
the same formula still satisfies the pentagon equation but it is
only a partial isometry. The purpose of the present paper is to
give necessary and sufficient conditions for an operator $V\colon
\Hil\o\Hil\to\Hil\o\Hil$ to determine a $C^*$-weak Hopf algebra.
 
$C^*$-weak Hopf algebras (WHA) are finite dimensional "quantum
groups"
with coproduct, counit, and antipode, but have no 1-dimensional
representations in general. Thus the counit is not an algebra map
and the antipode axioms have to be weakened accordingly. For its
axioms see \cite{BSz,Sz} and for a detailed exposition of these
quantum groups we refer to \cite{BNSz}. The main advantage of WHA's
in describing, for instance, the symmetry of the superselection
sectors in low dimensional QFT, is the flexibility of their
representation theory. Given any rigid monoidal $C^*$-category
{\bf C} with finitely many irreducible objects one can construct
a $C^*$-weak Hopf algebra $A$ with representation category
equivalent to {\bf C}. Roughly speaking this means that
$C^*$-WHA's exist for arbitrary (finite) set of 6j-symbols.
Since the 6j-symbols do not determine a unique $C^*$-WHA, one has
to supply more data than just a category. These data are provided
for example by a finite index depth 2 inclusion $N\subset M$ of
von Neumann algebras with finite dimensional centers
\cite{NSzW}. For II$_1$ factors and weak Kac algebras see
\cite{Nik,NV}.
 
In a recent paper \cite{EV} M. Enock and J.-M. Vallin study the
situation of
a general depth 2 inclusion of von Neumann algebras with a regular
operator valued weight and construct a certain isometry
called a pseudo-multiplicative unitary \cite{V}. In the finite
index case it is worth to compare their construction with ours. In
Section 6 we discuss the relation of finite dimensional
pseudo-multiplicative unitaries to multiplicative isometries and
reveal also some connection with Ocneanu's non-Abelian cohomology
\cite{O}. It will be shown that a unital multiplicative partial
isometry $V\colon\Hil\o\Hil\to\Hil\o\Hil$, what we introduce
in Sections 2 and 3, always determines a pseudo-multiplicative
unitary $U\colon\Hil\oR\Hil\to\Hil\oL\Hil$. By the results of
Section 3 this situation corresponds to the case when the `right
leg' and `left leg' of $V$, the algebras $A$ and $\duA$,
respectively, are weak bialgebras in the sense used in
\cite{BNSz}.
In Section 4 we put stronger conditions on $V$ and assume that it
satisfies a regularity condition, generalizing the one of
\cite{BaSk}.
Then we show that $A$ and $\duA$ are $C^*$-weak Hopf algebras in
duality.
 
The way from pseudo-multiplicative unitaries to multiplicative
isometries is not completely understood. Although we show at the
end of Section 6 that every $U$ determines a multiplicative
isometry $V$, unitalnes or regularity of this $V$ remain
unresolved.

\section{Multiplicative partial isometries}
Let $\Hil$ be a Hilbert space and $V\colon\Hil\o\Hil\to\Hil\o\Hil$
be a partial isometry, i.e. $VV^*V=V$. We shall say that $V$ is a
{\em multiplicative partial isometry} (MPI) if the following
equations hold on the 3-fold tensor product $\Hil\o\Hil\o\Hil$:
\bea
V_{23}V_{12}&=&V_{12}V_{13}V_{23}     \label{penta}\\
V_{13}V_{23}V_{23}^*&=&V_{12}^*V_{12}V_{13} \label{hexa1}\\
V_{12}V_{12}^*V_{23}&=&V_{23}V_{12}V_{12}^* \label{hexa2}\\
V_{12}V_{23}^*V_{23}&=&V_{23}^*V_{23}V_{12}\ .\label{hexa3}
\eea
The following equations are immediate consequences:
\bea
V_{12}^*V_{23}V_{12}&=&V_{13}V_{23} \label{first}\\
V_{23}V_{12}V_{23}^*&=&V_{12}V_{13} \label{second}\\
V_{12}V_{23}^*&=&V_{23}^*V_{12}V_{13} \label{third}\\
V_{12}^*V_{23}&=&V_{13}V_{23}V_{12}^* \label{fourth}\\
V_{12}V_{13}V_{13}^*&=&V_{23}V_{23}^*V_{12} \label{fifth}\\
V_{13}^*V_{13}V_{23}&=&V_{23}V_{12}^*V_{12}\ .\label{sixth}
\eea
For example, in order to obtain (\ref{first}) multiply
(\ref{penta}) by $V_{12}^*$ and then use (\ref{hexa1}). The reader
may easily prove the remaining equations in order of appearence.
For a geometrical interpretation of these equations see Section 6.
 
In this note we restrict ourselves to MPI's on finite dimensional
$\Hil$. Let $\LH$ denote the space of linear operators on
$\Hil$ and $\LH_*$ the space of linear functionals on $\LH$.
Let $V$ be any operator $V\in\LH\o\LH$ and construct the linear
maps
\bea
\lambda\colon\LH_*\to\LH\qquad &&\rho\colon\LH_*\to\LH\\
\lambda(\omega):=(\omega\o\id)(V)\qquad&&\rho(\omega)
:=(\id\o\omega)(V)
\eea
Their images $A:=\lambda(\LH_*)$ and $\duA:=\rho(\LH_*)$, called
the right leg and left leg of $V$, respectively,
are subspaces of $\LH$ that are in duality with respect to the
non-degenerate bilinear form
\beq
\bra\lambda(\omega),\rho(\omega')\ket\ :=\ (\omega\o\omega')(V)
\equiv\omega(\rho(\omega'))\equiv\omega'(\lambda(\omega))\ .
\label{pair}
\eeq
One obtains directly that $V\in\duA\o A$.
 
Let us introduce the following two binary operations on $\LH_*$.
\beq                  \label{convol}
\left.\begin{array}{l}
(\omega\star\omega')(X)\ :=\ (\omega\o\omega')(V^*(\one\o X)V)\\
(\omega\diamond\omega')(X)\ :=\ (\omega\o\omega')(V(X\o\one)V^*)
\end{array}\right\}\quad X\in\LH
\eeq
If $V$ is an MPI then we obtain
\bea
\lambda(\omega)\lambda(\omega')&=&(\omega\o\omega'\o\id)(V_{13}
V_{23})\eq{first}\lambda(\omega\star\omega')\label{A algebra}\\
\rho(\omega)\rho(\omega')&=&(\id\o\omega\o\omega')(V_{12}
V_{13})\eq{second}\rho(\omega\diamond\omega')\label{duA algebra}
\eea
showing that $A$ and $\duA$ are subalgebras of $\LH$.
 
The next step is to introduce the would-be coproducts $\cop$ and
$\ducop$, at first as linear maps $\LH\to\LH\o\LH$,
\beq\left.\begin{array}{l}
\cop(X)\ :=\ V(X\o\one)V^*\\
\ducop(X)\ :=\ V^*(\one\o X)V
\end{array}\right\}\quad X\in\LH\ .\label{deltas}\eeq
\begin{lem} \label{lem cop}
$\cop$ and $\ducop$ restrict to algebra maps $\cop\colon A\to A\o
A$ and $\ducop\colon\duA\to\duA\o\duA$.
\end{lem}
\Proof The identities
\beanon
\cop(\lambda(\omega))&=&(\omega\o\id\o\id)(V_{23}V_{12}V_{23}^*)
\eq{second}(\omega\o\id\o\id)(V_{12}V_{13})\ \in\ A\o A\\
\ducop(\rho(\omega))&=&(\id\o\id\o\omega)(V_{12}^*V_{23}V_{12})
\eq{first}(\id\o\id\o\omega)(V_{13}V_{23})\ \in\ \duA\o\duA
\eeanon
show that $\cop(A)\subset A\o A$ and
$\ducop(\duA)\subset\duA\o\duA$ so we have the required
restrictions. It remains to show multiplicativity of these
restrictions.
\beanon
&&\cop(\lambda(\omega))\cop(\lambda(\omega'))=
(\omega\o\omega'\o\id\o\id)(V_{13}V_{14}V_{23}V_{24})\eq{first}\\
&=&(\omega\o\omega'\o\id\o\id)(V_{12}^*V_{23}V_{12}V_{12}^*V_{24}V_{12})
\eq{hexa2}\\
&=&(\omega\o\omega'\o\id\o\id)(V_{12}^*V_{23}V_{24}V_{12})=
((\omega\star\omega')\o\id\o\id)(V_{12}V_{13})\eq{second}\\
&=&\cop(\lambda(\omega)\lambda(\omega'))
\eeanon
\beanon
&&\ducop(\rho(\omega))\ducop(\rho(\omega'))=
(\id\o\id\o\omega\o\omega')(V_{13}V_{23}V_{14}V_{24})\eq{second}\\
&=&(\id\o\id\o\omega\o\omega')(V_{34}V_{13}V_{34}^*V_{34}V_{23}V_{34}^*)
\eq{hexa3}\\
&=&(\id\o\id\o\omega\o\omega')(V_{34}V_{13}V_{23}V_{34}^*)=
(\id\o\id\o(\omega\diamond\omega'))(V_{13}V_{23})\eq{first}\\
&=&\ducop(\rho(\omega)\rho(\omega'))
\eeanon
\qed
 
From now on $\cop$ and $\ducop$ will denote these restrictions of
the original maps (\ref{deltas}).
\begin{lem} \label{lem dual}
Under the pairing $\bra\ ,\ \ket$ the comultiplication maps $\cop$
and $\ducop$ are the transposes of the multiplications on $\duA$
and $A$, respectively. In particular $\cop$ and $\ducop$ are
coassociative.
\end{lem}
\Proof We need to show that for $\omega,\omega',\omega''\in\LH_*$
\beanon
\bra\lambda(\omega),\rho(\omega')\rho(\omega'')\ket&=&
\bra\cop(\lambda(\omega)),\rho(\omega')\o\rho(\omega'')\ket\\
\bra\lambda(\omega)\lambda(\omega'),\rho(\omega'')\ket&=&
\bra\lambda(\omega)\o\lambda(\omega'),\ducop(\rho(\omega''))\ket
\eeanon
or, equivalently
\beanon
(\omega'\diamond\omega'')(\lambda(\omega))&=&(\omega\o\omega'\o\omega'')
(V_{12}V_{13})\\
(\omega\star\omega')(\rho(\omega''))&=&(\omega\o\omega'\o\omega'')
(V_{13}V_{23})
\eeanon
which, up to an application of (\ref{second}) or (\ref{first}),
are precisely the definitions of the convolution products
(\ref{convol}).\qed
 
In this way we have shown that a multiplicative partial isometry
determines
a pair $(A,\duA)$ of algebras in duality such that the induced
comultiplications are algebra maps. It is not clear, however, if
these algebras have units or if they are closed under the
$^*$-operation. So we need further assumptions.
 

\section{Unital MPI's and Weak Bialgebras}
 
At first we will seek for the conditions on the finite
dimensional MPI $V$ that ensure that $A$ and $\duA$ are weak
bialgebras (WBA's) in the sense of \cite{BNSz}. Obviously it is
necessary that both of them should be unital algebras
(hence counital coalgebras). We claim that this condition, called
{\em unitalness}, is not only necessary but also sufficient. It
is also shown that under this condition the elements of $A$ and
$\duA$ realize a (not necessarily faithful) representation of the
Weyl algebra (or Heisenberg double) $A\cros\duA$
\cite{BSz,BNSz}.
 
\begin{defi} \label{unital}
A finite dimensional MPI $V$ on the Hilbert space $\Hil$
is {\em unital} if there
exist functionals $\LH_*\ni\eps$  and $\dueps$ such that
$A\ni \lambda(\dueps)\equiv \1$ and $\duA \ni \rho (\eps) \equiv \du1$
are two-sided units for $A$ and $\duA$, respectively.
\end{defi}
 
In order to illustrate that, in contrast to multiplicative
unitaries, finite dimensional MPI's are not always unital, let
stand here a non-unital example.
Let $\Hil=\C^2$ and define $V=e_{11}\o e_{12}+
e_{22}\o e_{22}$ with a chosen set of $^*$-matrix units
$\{e_{ij}\}_{i,j\in\{1,2\}}$. Then one can see by inspection that
$V$ is an MPI, its left leg contains $\1$, but its right leg does
not.
 
Although the functionals $\eps$ and $\dueps$ in the above
Definition are not unique they have a unique restriction onto $A$
and $\duA$, respectively. These restrictions (also denoted as
$\eps$ and $\dueps$) are then counits of $A$ and $\duA$,
respectively.
 
If $V$ is unital then $A$ and $\duA$ are WBA's provided the counits
are weakly multiplicative or, equivalently,
if the units $\du1$ and $\1$ are weakly comultiplicative. We show
this latter property using
 
\begin{lem} \label{delt1}
Let $V$ be a finite dimensional unital MPI on the Hilbert space $\Hil$
with unit elements $\1\in A$ and $\du1\in\duA$. Then
\bea  \cop (\1)&=& VV^* \label{cop1}\\
      \ducop (\du1) &=& V^*V.\label{ducop1}
\eea
\end{lem}
\Proof By (\ref{deltas}) we have for any $\omega\in \LH_*$
\bea (\one\o\lambda(\omega))\cop(\1)&=& (\omega\o\id\o\id)(V_{13}V_{23}
(\one\o\1\o\one ) V_{23}^* )\eq{first}\nn
&=& (\omega\o\id\o\id)(V_{12}^* V_{23} V_{12} (\one\o\1\o\one)
V_{23}^*) . \nonumber \eea
By  the assumption that $\1$ is a (right) unit for $A$,
$V(\one\o\1)=V$ and
\[
(\one \o \lambda(\omega))\cop (\1)=(\omega\o\id\o\id)(V_{13} V_{23}
V_{23}^*) =(\one \o\lambda(\omega)) VV^* .
\]
Setting $\omega=\dueps$ and using the assumption that $\1$ is a (left)
unit for $A$ (\ref{cop1}) is proven.
A similar argument shows that
\beq  \ducop (\du1) (\rho(\omega)  \o \one)= V^*V (\rho(\omega)  \o \one)
\nonumber\eeq
for any $\omega\in\LH_*$, hence the substitution $\omega=\eps$
proves (\ref{ducop1}). \qed
 
As a consequence of Lemma \ref{delt1}
\beq (\cop(\1)\o \1)(\1\o \cop (\1))=V_{12}V_{12}^* V_{23}V_{23}^*
\eq{hexa2} V_{23} V_{23}^* V_{12} V_{12}^* =
(\1\o \cop (\1))(\cop(\1)\o \1) \nonumber \eeq
which, by (\ref{fifth}), equals to
\beq \1\c\o\1\cc\o\1\ccc = V_{12} V_{13} V_{13}^* V_{12}^* .\nonumber\eeq
Similarly,
\beq (\ducop(\du1)\o\du1)(\du1\o\ducop(\du1))=V_{12}^* V_{12} V_{23}^* V_{23}
\eq{hexa3} V_{23}^* V_{23} V_{12}^* V_{12}= (\du1\o \ducop(\du1))
(\ducop(\du1)\o \du1), \nonumber\eeq
which, by (\ref{sixth}), equals to
\beq \du1\c\o\du1\cc\o\du1\ccc = V_{23}^* V_{13}^* V_{13} V_{23}. \nonumber\eeq
This proves that if $V$ is unital then the resulting algebras $A$ and $\duA$
are WBA's in duality.
 
A further consequence of the above Lemma is that the subalgebras
$A^L$ and $A^R$ of $A$, that were originally defined as the right
leg and left leg, respectively of $\cop(\1)$ \cite{BNSz}, appear
now in the form
\bea
A^L&=&\{\,(\omega\o\id)(VV^*)\,|\,\omega\in\LH_*\,\}\\
A^R&=&\{\,(\id\o\omega)(VV^*)\,|\,\omega\in\LH_*\,\}\ .
\eea
Therefore they are selfadjoint subalgebras in $\LH$ even if we do
not know whether $A$ is selfadjoint. Similar conclusion holds for
the subalgebras $\duA^L$ and $\duA^R$ of $\duA$.
 
As far as the relative position of $A$ and $\duA$ in $\LH$ is
concerned we want to show that $A$ and $\duA$
generate a representation of the Weyl algebra $A\cros \duA$ on $\Hil$.
As a matter of fact the
pentagon equation (\ref{penta}) implies the commutation relation
\bea \label{Weyl rel}
\rho(\omega)\lambda(\omega')&=&(\omega'\o\id\o\omega)(V_{23}V_{12})=\nn
  &=&(\omega'\o\id\o\omega)(V_{12}V_{13}V_{23})=
  (\id\o\omega)(\cop(\lambda(\omega'))V)=\nn
  &=&\lambda(\omega')\c\bra\lambda(\omega')\cc,\rho(\omega)\c\ket
  \rho(\omega)\cc.
\eea
The only missing Weyl algebra relation is $\1=\du1$.
\begin{prop} \label{weyl}
Let $V$ be a finite dimensional unital MPI on the Hilbert space $\Hil$
with unit elements $\1\equiv \lambda(\dueps)\in A$
and $\du1\equiv \rho(\eps)\in \duA$. Then
\beq \1=\du1 \label{coinc}\eeq
as elements of $\LH$.
\end{prop}
\Proof We recall \cite{BNSz} that a projection from $\duA$ onto
$\duA^L$
is provided by $\duPL(\varphi):=\dueps(\du1\c\varphi)\du1\cc$.
Hence for arbitrary $\omega\in\LH_*$
\bea \duPL (\rho(\omega))&=&
(\dueps \o \id\o \omega)(V_{12}^* V_{12} V_{13}) \eq{hexa1}
(\dueps\o\id\o\omega )(V_{13} V_{23} V_{23}^*)=\nn
&=& (\id\o \omega)((\one\o\1)VV^*)=(\id\o\omega)(\cop(\1)).
\label{am}\eea
Setting $\omega=\eps$ we obtain (\ref{coinc}). \qed
 
As a byproduct equation (\ref{am}) tells us that the
subalgebras $A^R\subset A$ and $\duA^L\subset\duA$ coincide as
subalgebras of $A\cros\duA$ and therefore of $\LH$. As a
counterpart of this relation one can also show that
\beq                        \label{am2}
\PR(\lambda(\omega))\equiv\1\c\eps(\lambda(\omega)\1\cc)=
(\omega\o\id)(V^*V)=\bra\lambda(\omega),\du1\c\ket\du1\cc\ ,
\eeq
hence $A^R=\duA^L$ and the identification is given by $A^R\ni
x^R\mapsto (\du1\la x^R)\in\duA^L$. This relation is called the
amalgamation relation.

\section{Regular MPI's and Weak Hopf Algebras}
 
Given a finite dimensional unital MPI $V$ and the associated WBA's
$A$ and $\duA$ one may look for the extra conditions on $V$ that
ensure one of the following special cases to occur:
\begin{itemize}
\item There exist antipodes $S$ and $\duS$ making $A$ and $\duA$
weak Hopf algebras.
\item $A$ and $\duA$ are closed under the $*$-operation.
\item $A$ and $\duA$ are $C^*$-WHA's in duality.
\end{itemize}
It turns out that these cases occur at the same time. In this
Section we give a necessary and sufficient condition for this to
happen that is reminiscent to the regularity condition of
\cite{BaSk}.
 
{\em Remark}\,: Questions like whether $A$ and $\duA$ are
selfadjoint do not occur in the works \cite{BaSk,EV}.
In their approach the Hopf algebra (Hopf bimodule) is {\em
defined} to be the selfadjoint closure of the right or left
leg of the (pseudo-) multiplicative unitary. In our finite
dimensional approach the WBA or WHA $A$ is the right leg of the
MPI $V$ and not larger. On the one hand this is very natural
in view of the duality of $A$ and $\duA$ under the pairing
(\ref{pair}) but on the other hand this will cause difficulties
if one wants to compare MPI's with pseudo-multiplicative unitaries
(see Section 6).
 
\begin{prop} Let $V$ be a finite dimensional MPI on the Hilbert space
$\Hil$ such that the resulting
algebras $A$ and $\duA$ are WHA's
with coproducts given in (\ref{deltas})
and with (the unique) antipodes $S\colon A\to  A$ and
$\duS:\duA\to \duA$. Then we have the relation
\beq V^*= (\duS\o\id)(V)\equiv(\id\o S)(V) \eeq
and therefore $A$ and $\duA$ are $^*$-subalgebras of $\LH$.
\end{prop}
 
\Proof (\ref{pair}) implies that $V=\sum_i \beta^i \o b_i$ with any basis
$\{b_i\}$ of $A$ and its dual basis $\{\beta^i\}$ of $\duA$. Let
$V\p\colon =\sum_i \duS(\beta^i)\o b_i\equiv \sum_i \beta^i\o S(b_i)$.
We claim that $V\p=V^*$. Using the assumption that $A$ and $\duA$ are
WHA's in duality compute
\bea VV\p &=& \sum_{i,j} \beta^i \duS(\beta^j) \o b_i b_j=
\sum_k \duPL(\beta^k)\o b_k =\nn
&=& \du1\cc \o \1\ra\du1\c = \du1\ra\1\c\o\1\cc=\1\c\o\1\cc=VV^*\nn
V\p V&=& \sum_{i,j} \duS(\beta^i)\beta^j \o b_ib_j =
\sum_k \duPR(\beta^k) \o b_k =\nn
&=& \du1\c\o\du1\cc\la\1= \du1\c\o\du1\cc=V^*V
\nonumber\eea
where in the last step of both cases we used the amalgamation
relation (\ref{am}-\ref{am2}). Now
\beanon
V^*&=&V^*VV^*=V\p VV^*=
V\p \cop(\1)=\sum_i \beta^i\1\c\o S(b_i)\1\cc=\nn
&=& \sum_i \beta^i\ra\1\c\o S(b_i)\1\cc=
\beta^i\o S(b_i) S(\1\c)\1\cc = V\p
\eeanon
therefore
\beq\left.
\begin{array}{l}
 S(\lambda(\omega))\ =\ (\omega\o\id)(V^*)\\
\duS(\rho(\omega))\ =\ (\id\o\omega)(V^*)
\end{array}
\right.\label{antipodes}\eeq
implying that $S(A)\subset A^*$ and $\duS(\duA)\subset\duA^*$.
This is possible for the bijections $S:A\to A$ and
$\duS:\duA\to\duA$
only if $A$ and $\duA$ are $^*$-subalgebras of $\LH$. \qed
 
The next Proposition proves a converse result plus some more.
\begin{prop} \label{ap}
Suppose that the MPI $V$ on the Hilbert space $\Hil$
is such that its right and left leg, $A$ and $\duA$,
are $^*$-subalgebras of $\LH$. Then $V$ is unital and the
expressions (\ref{antipodes}) define antipodes that
make $A$ and $\duA$ $C^*$-WHA's in duality.
\end{prop}
 
\Proof Since $^*$-subalgebras of $\LH$ are semisimple, $A$ and
$\duA$ have units. Furthermore, being in duality by the pairing
(\ref{pair}), they possess functionals $\eps$ and $\dueps$
required in Definition \ref{unital}. Thus $V$ is unital and $A$
and $\duA$ are WBA's in duality by the results of Section 3.
 
In order to construct antipodes notice that if $\lambda(\omega)=0$
then $\omega(\duA^*)=\omega(\duA)=0$, therefore the $S$ of
(\ref{antipodes}) is a well defined map $A\to A$. Similarly,
(\ref{antipodes}) defines a map $\duS\colon\duA\to\duA$. These
maps are the transpose of each other with respect to the canonical
pairing (\ref{pair}),
\[
\bra\duS(\rho(\omega)), \lambda(\omega\p)\ket =
(\omega\p\o \omega)(V^*)=
\bra\rho(\omega), S(\lambda(\omega\p))\ket
\]
for all $\omega,\omega\p\in\LH_*$. It remained to show that the
$C^*$-WHA axioms are satisfied.
 
Define the antilinear involution $_*\colon\LH_*\to\LH_*$ by
$\omega_*(X)\ :=\ \overline{\omega(X^*)}$ for $X\in\LH$.
Then (\ref{antipodes}) can be rewritten as
\beanon
S(\lambda(\omega))&=&\lambda(\omega_*)^* \\
\duS(\rho(\omega)) &=&\rho(\omega_*)^*\ .
\eeanon
By showing that $_*$ preserves both convolution products
(\ref{convol}),
\[
(\omega\star\omega')_*=\omega_*\star\omega'_*\ ,\quad
(\omega\diamond\omega')_*=\omega_*\diamond\omega'_*\ ,
\]
we find that both $S$ and $\duS$ are anti-multiplicative and
anti-comultiplicative. Finally
\beanon
\lambda(\omega)\c\o\lambda(\omega)\cc S(\lambda(\omega)\ccc)&=&
(\omega\o\id\o\id)(V_{12} V_{13} V_{13}^*)\eq{fifth}\\
&=&(\omega\o\id\o\id)(V_{23}V_{23}^*V_{12})=\cop(\1)
(\lambda(\omega)\o\one)
\\
\rho(\omega)\c \o \rho(\omega)\cc \duS (\rho(\omega)\ccc) &=&
(\id \o \id \o \omega) (V_{13} V_{23} V_{23}^*)\eq{hexa1}\\
&=&(\id\o\id\o\omega)(V_{12}^* V_{12} V_{13})={\hat\Delta}(\du1)
(\rho(\omega)\o\one)
\eeanon
prove that the WHA axioms of \cite{Sz} hold both in $A$ and
$\duA$. Since the coproducts (\ref{deltas}) are manifestly $^*$-
algebra maps, $A$ and $\duA$ are $^*$-WHA's. Furthermore
the defining representations of $A$ and $\duA$ on $\Hil$ are
faithful $^*$-representations by construction therefore $A$ and
$\duA$ are $C^*$-WHA's. \qed
 
It remains to characterize the situation of $A$ and $\duA$ being
selfadjoint in "more algebraic" terms, i.e. using only the
relative positions of $A$ and $\duA$ in $\LH$ without referring to
their $^*$-structure. This will
be the regularity condition on the multiplicative isometry $V$.
 
In analogy with \cite{BaSk} we define the subspace
${\cal C}(V):= V_2 \LH V_1$ in $\LH$, where $V_1\o V_2$ stands
for $V$, and verify using the pentagon equation (\ref{penta})
that ${\cal C}(V)$ is a subalgebra of $\LH$.
 
\begin{lem} \label{c*}
Let $V$ be a unital MPI on the Hilbert space $\Hil$.
If ${\cal C}(V)$ is a $^*$-subalgebra of $\LH$ then so are $A$
and $\duA$.
\end{lem}
 
\Proof  The proof generalizes the one of Proposition 3.5 in \cite{BaSk}.
At first we show that
\beq A^*=\{ (\omega\o\omega\p\o \id)(\Sigma_{12} V^*_{23} V_{12} V_{13})
\vert \omega,\omega\p\in\LH_*\}
\label{a*} \eeq
where $\Sigma:\Hil\o\Hil\to\Hil\o\Hil$ is the flip map.
This follows from the computations
\beanon
(\omega\o\omega\p\o\id)(\Sigma_{12} V_{23}^* V_{12} V_{13})
&\eq{third}& (\omega\o\omega\p\o\id)(\Sigma_{12} V_{12} V_{23}^*)=\\
&=& (\omega\p\o\id )( ((\omega\o\id)(\Sigma V)\o\one) V^* ) \in A^*
\eeanon
and
\[
(\omega\o \id)(V^*)= (\omega_1 \o \omega_2\o\id)
(\Sigma_{12} V^*_{23} V_{12} V_{13})
\]
where we introduced $\omega_1\o\omega_2\in\LH_*\o\LH_*$ by setting
$(\omega_1\o\omega_2)(X):=(\dueps \o\omega)(\Sigma X)$.
Thus (\ref{a*}) is proven.
The next identity
\beanon
(\omega\o\omega\p\o\id)(\Sigma_{12} V_{23}^* V_{12} V_{13}) &=&
 (\omega\o\omega\p\o\id)(V_{13}^* \Sigma_{12} V_{12} V_{13}) =\\
&=&(\omega\o\id)(V^*((\id\o\omega\p)(\Sigma V)\o\one)V)
\eeanon
shows that if ${\cal C}(V)\equiv \{(id\o \omega)(\Sigma V)\vert
\omega\in \LH_*\}$
is closed under the $^*$-operation then so is $A^*$ hence $A$.
 
In the case of $\duA$  repeat the above argument using the
fact that
in passing from the MPI $V$ to the MPI $\Sigma V^*\Sigma$ the left
leg $\duA (V)$ becomes the adjoint of the
right leg $A(\Sigma V^*\Sigma)$ and also
${\cal C}(\Sigma V^* \Sigma) ={\cal C}(V)^*$. \qed
 
\smallskip
$A^R$ is the subalgebra of $\LH$ spanned by the elements
$\{(\omega\o\id)(V^*V)\,|\,\omega\in\LH_*\,\}$. It is
obviously a $^*$-subalgebra and for $a^R=(\omega\o\id)(V^*V)$
\beanon
(\one\o a^R)V&=&(\omega\o\id\o\id)
(V_{13}^* V_{13} V_{23}) \eq{sixth}\\
&=& (\omega\o\id\o\id)(V_{23} V_{12}^* V_{12})=
V (a^R\o\one)
\eeanon
hence $A^R$ commutes with ${\cal C}(V)$. Let us make
the following
\begin{defi} \label{regular}
A finite dimensional unital MPI $V$ on the Hilbert space $\Hil$ is
called {\em regular} if
\beq {\cal C}(V)= (A^R)\p \cap \1\LH\1\ . \eeq
\end{defi}
 
In the special case of $V$ being a multiplicative unitary the
$A^R$ consists only of the scalars therefore $(A^R)\p\cap\1\LH\1
=\LH$ and our regularity condition reduces to the regularity of
\cite{BaSk}. Although in finite dimensions all multiplicative
unitaries are regular by Theorem 4.10 of \cite{BaSk} we do not
know any generalization of this result to multiplicative
isometries.
 
\begin{thm} \label{thm}
The algebras $A$ and $\duA$ obtained from a finite dimensional MPI
$(V,\Hil)$ are $^*$-subalgebras of $\LH$
if and only if $V$ is unital and regular.
\end{thm}
 
\Proof Since $A^R$ is a $^*$-subalgebra of $\LH$ so is its commutant.
This implies that if $V$ is unital and regular then ${\cal C}(V)$ is
a $^*$-subalgebra of
$\LH$ and using Lemma \ref{c*} the {\em if} part follows.
 
To prove the converse statement suppose that $A$ and $\duA$ are
$^*$-subalgebras
of $\LH$ so they are $C^*$-WHA's in duality by Proposition \ref{ap}.
Then $V$ is necessarily unital. In this case
$A^R$ is the right subalgebra of $A$ coinciding with the left subalgebra
of $\duA$ (see (\ref{am})).
 
Knowing already that ${\cal C}(V)\subset (A^R)\p \cap \1\LH\1$ it
remains to show that
also $(A^R)\p \cap \1\LH\1\subset {\cal C}(V)$. For that purpose
let $X\in (A^R)\p \cap \1\LH\1$. Then
\beanon
X&=& X\1= X\1\c S(\1\cc)=X\1\c \duS\inv (\1\cc\la \du1)=\\
&=& X \sum_k \PR(b_k) \duS\inv (\beta^k)= \sum_k
\PR(b_k)X\duS\inv(\beta^k)=\\
&=&\sum_kS(b_{k(1)})b_{k(2)}X \duS\inv(\beta^k)=
 \sum_{i,j} S(b_i) b_j X \duS\inv(\beta^i \beta^j)=\\
&=&\sum_{i,j} b_i b_j X \duS\inv (\beta^j) \beta^i \in {\cal C}(V)
\eeanon
finishes the proof.\qed
 
\smallskip
With the above Theorem
we have characterized the class of MPI's that lead to $C^*$-WHA's.
The question arises whether all $C^*$-WHA's can be obtained in this
way. The answer is in fact very easy.
Let $A$ be a $C^*$-weak Hopf algebra and let $\pi\colon
A\cros\duA\to\LH$ be a $^*$-representation such that the
restrictions $\pi|_A$ and $\pi|_{\duA}$ are faithful. Choose
a basis $\{b_i\}$ of $A$ and construct the dual basis
$\{\beta^i\}$, $\bra b_i,\beta^j\ket=\delta_{ij}$ of $\duA$. Then
\beq
V\ :=\ \sum_i\ \pi(\beta^i)\,\o\,\pi(b_i)
\eeq
is a multiplicative partial isometry in the sense of
(\ref{penta}), (\ref{hexa1}), (\ref{hexa2}), and
(\ref{hexa3}) and furthermore it is unital and regular in the
sense of Definitions \ref{unital} and \ref{regular}. The proof
of this statement is
an elementary weak Hopf calculus which we omit. Notice
that as a special case we obtain the "classical" example when
$\Hil$ is the left regular representation of a $C^*$-WHA $A$ with
scalar product provided by the Haar measure, $(x,y)=\bra x^*y,\hat
h\ket$, $x,y\in A$. In this case the action of $V$ is given by
$V(x\o y)=x\c\o x\cc y$.

\section{Pseudo-multiplicative unitaries in finite dimensions}
 
In order to discuss the relation of MPI's to the
pseudo-multiplicative unitaries \cite{V,EV}
we specialize their definition to the case when the
Hilbert space in the game is finite dimensional. At first we
exhibit the Connes-Sauvageot relative tensor product \cite{S} of
finite
dimensional modules as a subspace in the ordinary tensor product.
Then the pseudo-multiplicative unitary $U$ will be obtained by
restricting the domain and range of the MPI $V$ to its initial and
final support. It should be emphasized, however, that
the pseudo-multiplicative unitary has to be supplied with
an a priori knowledge of the algebra $A^L$ and a faithful state
on it while this information is implicitely stored in the
structure of $V$.
 
Let $B$ be a finite dimensional $C^*$-algebra, ${\cal H}$ and
${\cal K}$ finite dimensional Hilbert spaces, ${\cal H}$ carrying
a right and ${\cal K}$ a left $B$-module structure, i.e. there
are given $^*$-homomorphisms $\beta:B^{o}\to \LH$ and
$\gamma:B\to
{\cal L}({\cal K})$. If $\psi:B\to \C$ is a faithful positive
linear functional the {\em relative tensor product}
of $\Hil_{\beta}$ and $_{\gamma}{\cal K}$ over
$\psi$ is defined to be the subspace in $\Hil\o{\cal K}$
obtained as the
image of a projection $E_\psi\in{\cal L}(\Hil\o{\cal K})$
constructed below.
 
Let $\{a_i\}$ be a basis of $B$ and $\{b_i\}$ the dual basis
with respect to $\psi$, i.e. $\psi(b_i a_j)=\delta_{i,j}$. Then
$x=\sum_i a_i \psi(b_ix)$ for all $x\in B$, i.e.
$\{a_i,\,b_i\}$ is a {\em quasibasis} of $\psi$ in the sense of
\cite{Wata}. The {\em index} of $\psi$,
$\lambda:=\sum_i a_i b_i$, is a positive
invertible element of $\Center B$. The {\em modular
automorphism} of $\psi$ is the (non-$^*$) automorphism
$\theta_{\psi}$ of $B$ satisfying $\psi(xy)=
\psi(y\theta_{\psi}(x))$ for all $x,y\in B$.
In terms of these data we can define an element $e_\psi\in
B^o\o B$ by the formula
\beq                                        \label{e}
e_\psi\equiv\sum_i u_i\o v_i:=
\sum_i \lambda\inv a_i\o\theta_\psi^{1/2}(b_i)\ .
\eeq
Checking that $e_\psi$ is a Hermitean idempotent we have
$E_\psi:=(\beta\o\gamma)(e_\psi)$ as the projection defining the
relative tensor product $\Hil\amalgo{\psi}{\cal
K}:=E_\psi(\Hil\o{\cal K})$.
The image of $\xi\o\eta\in\Hil\o{\cal K}$ in the relative tensor
product will be denoted by $\xi\amalgo{\psi}\eta$. Using the
property $\sum_i a_i\o b_ix=\sum_i xa_i\o b_i$,
$x\in B$, of the quasibasis
we immediately obtain the amalgamation relation
\beq \label{heart1}
\xi\amalgo{\psi}\gamma(x)\eta=\beta\circ\theta_{\psi}^{-1/2}(x)\xi
\amalgo{\psi}\eta
\eeq
for all $\xi\amalgo{\psi}\eta\in\HK$.
 
The above definition of the relative tensor product applies also to $\KH$
if we replace $B$ with $B^o$ and call the resulting functional
$\psi^o$. The identities $\sum_i b_i\psi(xa_i)=x$ and $\psi(yx)=
\psi(\theta_{\psi}{\inv}(x) y)$ show that $\sum_i a_i^o\o b_i^o\colon =
\sum_i b_i \o a_i$ is the quasi-basis of $\psi^o$ and
$\theta_{\psi^o}=\theta_{\psi}\inv$. Therefore $\sum_i u_i^o\o v_i^o =
\sum_i u_i\o v_i$ and $\KH$ is defined by the projection
$E_{\psi^o}=(\gamma\o\beta)(e_\psi)$.
Denoting the image of $\eta\o\xi$ in $\KH$ by $\eta\amalgo{\psi^o}\xi$,
we obtain the amalgamation
\beq
\eta \amalgo{\psi^o} \beta(x)\xi = \gamma\circ \theta_{\psi}^{\half}(x)\eta
\amalgo{\psi^o} \xi
\label{heart2} \eeq
for all $\eta\amalgo{\psi^o}\xi \in \KH$.
 
Some caution is in order with the equations (\ref{heart1}) and
(\ref{heart2}).
They must not be considered as `the operator $\one\o\beta(x)$',
\dots etc, acting on $\eta\amalgo{\psi^o}\xi$. Rather the vectors
$\eta\o\beta(x)\xi\in {\cal K}\o \Hil$,\dots etc, are mapped into
the subspace $\KH$. Only operators $X\in{\cal L}({\cal K})$ commuting
with $\gamma(B)$ and $Y\in\LH$ commuting
with $\beta(B)$ can be restricted to operators
$(\one_{\Hil}\amalgo{\psi}X)$,
$(Y\amalgo{\psi}\one_{\cal K})\in {\cal L}(\HK)$ and
$(\one_{\cal K}\amalgo{\psi^o} Y),\ (X\amalgo{\psi^o} \one_{\Hil})\in
{\cal L}(\KH)$.
 
For later convenience we supress the letters $\beta$ and $\gamma$ and
write $\xi\cdot b$ and $b\cdot\eta$ for $\beta(b)\xi$ and $\gamma(b)\eta$,
respectively. In this spirit we may think $\amalgo{\psi}$ as the
symbol $\cdot u_i\o v_i\cdot$ (with the $i$ summed over).
 
The usual flip operator $\Sigma\colon\Hil\o{\cal K}\to{\cal
K}\o\Hil$
determines an isomorphism $\Sigma_\psi\colon\HK\to\KH$ by
restriction since $\Sigma$ intertwines between the projections
$E_\psi$ and $E_{\psi^o}$, or in other words, because $\sum_i
u_i\o v_i=\sum_i v_i\o u_i$. This follows using the fact that
the modular automorphism is necessarily
inner on a finite dimensional $C^*$-algebra. As a matter of
fact let $g_{\psi}\in B$ be a positive element implementing
$\theta_{\psi}$, i.e. $g_\psi x g_\psi\inv=\theta_\psi(x)$ for all
$x\in B$. Then
\beanon
\sum_i v_i\o u_i &=&
\lambda\inv \sum_i g_{\psi}^{\half}b_i g_{\psi}^{-\half}\o a_i=
\lambda\inv \sum_i g_{\psi}^{-\half} \theta_{\psi}(b_i) g_{\psi}^{\half}
\o a_i =\\
&=& \lambda\inv \sum_i g_{\psi}^{-\half} a_i g_{\psi}^{\half}\o b_i=
\sum_i u_i\o v_i.
\eeanon
 
With the above method one can construct also {\em multiple} relative
tensor products of modules over (different) finite dimensional
$C^*$-algebras.
Let $A$ and $B$ be finite dimensional $C^*$-algebras, $\Hil, {\cal K}$
and ${\cal M}$ Hilbert spaces with the following module
structures: $\cal H$ be a right $A$-module, ${\cal M}$ an
$A$-$B$ bimodule, and ${\cal K}$ a left $B$-module. Let $\phi:A\to\C$
and $\psi:B\to\C$ be faithful positive linear functionals. Then
there are two threefold relative tensor products defined
respectively by the formulae
\beanon
\Hil\amalgo{\phi}({\cal M}\amalgo{\psi} {\cal K})\colon =
\sum_{i,j} \Hil\cdot u^{\phi}_i \o v^{\phi}_i\cdot({\cal M}\cdot u^{\psi}_j
\o v^{\psi}_j\cdot {\cal K}) \\
(\Hil\amalgo{\phi} {\cal M})\amalgo{\psi} {\cal K} \colon =
\sum_{i,j} (\Hil\cdot u^{\phi}_i \o v^{\phi}_i\cdot {\cal M})\cdot u^{\psi}_j
\o v^{\psi}_j \cdot {\cal K},
\eeanon
which, as subspaces of $\Hil\o{\cal M}\o{\cal K}$, coincide up to
the associativity natural isomorphism in the category of Hilbert
spaces. Supressing this natural isomorphism we can denote this
Hilbert space by $\Hil\amalgo{\phi}{\cal M}\amalgo{\psi}{\cal K}$.
 
Considering the $A$-actions on $\Hil$ and ${\cal M}$ as $A^o$-actions
we have also the Hilbert space  $({\cal M}\amalgo{\psi}{\cal K})
\amalgo{\phi^o}\Hil$. Similarly one can define ${\cal
K}\amalgo{\psi^o}(\Hil\amalgo{\phi}{\cal M})$ and \break
${\cal K}\amalgo{\psi^o}{\cal M}\amalgo{\phi^o}\Hil$. They are all
naturally isomorphic under the flip maps:
 
\begin{picture}(250,180)(-40,0)
\put(0,150){$\Hil\amalgo{\phi}{\cal M}\amalgo{\psi}{\cal K}$}
\put(100,155){\vector(1,0){60}}
\put(120,165){$\Sigma_\phi$}
\put(200,150){$({\cal M}\amalgo{\psi}{\cal
                                      K}\,)\amalgo{\phi^o}\Hil$}
\put(25,130){\vector(0,-1){50}}
\put(5,105){$\Sigma_\psi$}
\put(0,50){${\cal K}\amalgo{\psi^o}(\Hil\amalgo{\phi}{\cal M})$}
\put(235,130){\vector(0,-1){50}}
\put(245,105){$\Sigma_\psi\amalgo{\phi^o}\one_{\Hil}$}
\put(100,55){\vector(1,0){60}}
\put(110,35){$\one_{\cal K}\amalgo{\psi^o}\Sigma_{\phi}$}
\put(200,50){${\cal K}\amalgo{\psi^o}{\cal
                                         M}\amalgo{\phi^o}\Hil$}
\end{picture}

\section{The relation of $U$ and $V$}
 
In this Section we will present two constructions. At first
we show how a finite
dimensional unital multiplicative isometry $(V,\Hil)$ determines a
pseudo-multiplicative unitary $U$. After that starting from a
finite dimensional pseudo-multiplicative unitary $U$ we construct
a MPI $V$.
 
Let $V$ be a unital MPI on the finite dimensional Hilbert
space $\Hil$, and $A$, $\duA$ the associated WBA's in duality,
both acting on $\Hil$. By Lemma \ref{delt1} the left and right
subalgebras of $A$ and of $\duA$ are selfadjoint subalgebras of
$\LH$. In particular $A^L$ is a $C^*$-algebra and the counit
$\eps$ restricts to a faithful positive functional on $A^L$.
 
We need the following facts from the theory of weak bialgebras
\cite{Nill, BNSz}. Although an antipode may not exist on $A$ we
can define a would-be-antipode $S$ on the subalgebra $A^LA^R$ by
setting $S(x^Lx^R):=\PL(x^R)\PR(x^L)\equiv(\du1\ra
x^L)\la\1\ra(x^R\la\du1)$. Then the element
$S(\1\c)\o\1\cc=\1\cc\o S^{-1}(\1\c)\in A^L\o A^L$ provides a
quasibasis of $\eps\colon A^L\to \C$, hence $\eps|_{A^L}$ has
index $\1$. The modular automorphism of $\eps|_{A^L}$ is
$\theta=S^2|_{A^L}$ and it is implemented by a positive element
$g_L\in A^L$. Although $g_L$ is not unique, the formulae
$\theta^{1/2}(x^L)=g_L^{1/2}x^Lg_L^{-1/2}$ and
$S_o:=S\circ\theta^{-1/2}$ do not depend on this ambiguity. Here
$S_o$ is a "unitary antipode" satisfying $S_o\circ\ ^*\,=\ ^*\circ
S_o$ and $S_o^2=\id$. By means of these definitions we can
construct
\beq
e_\eps:=\1\cc\o\theta^{1/2}(S^{-1}(\1\c))
\eeq
which is precisely the Hermitean idempotent (\ref{e}) needed in
relative tensor products of $A^L$-modules over $\eps$ or $\eps^o$.
 
Corresponding to the three $C^*$-subalgebras $A^L$,
$A^R\equiv\duA^L$, and $\duA^R$ of $\LH$ there are three mutually
commuting actions of $A^L$ on $\Hil$:
\beq \alpha_{01}(x^L)\xi\colon =x^L \xi \qquad
     \alpha_{02}(x^L)\xi\colon =S_o(x^L)\xi \qquad
     \alpha_{12}(x^L)\xi\colon =(x^L\la \du1) \xi
\eeq
$\alpha_{01}$ and $\alpha_{12}$ are left actions while
$\alpha_{02}$ is a right action. It is tempting to visualize this
trimodule structure of $\Hil$ by drawing a triangle $(012)$
for the Hilbert space
\beq                                       \label{trimod}
\Hil\ =\ \trimodule
\eeq
and say that the edge $(ij)$ is a left or right action of $A^L$
according to whether the relative orientation of $(ij)$ to the
2-simplex $(012)$ is positive or negative.
 
Now we want to exhibit the source and target spaces of the partial
isometry $V$ as relative tensor products of $\Hil$ with itself.
For that purpose we compute
\beanon
V^*V &=& \ducop(\du1)=\du1\c\o \du1\cc\la \1 =\1\cc\la\du1\o\1\c=\\
&=& \1\cc\la\du1 \o S_o(g_L^{\half}S\inv(\1\c)g_L^{-\half})=
(\alpha_{12}\o \alpha_{02})(e_\eps)\ ,\\
VV^* &=& \Delta(\1)= S_o(g_L^{\half} S\inv(\1\c) g_L^{-\half})\o\1\cc=\\
&=&(\alpha_{02}\o\alpha_{01})(e_\eps)\ .
\eeanon
This means that we may identify the source
and the target spaces of $V$ with the following relative tensor products:
\bea V^*V(\Hil\o\Hil) &=&
\source\label{prs}\\
     VV^*(\Hil\o\Hil) &=&
\target\label{prt}.
\eea
As a graphical representation of these relative tensor products
one draws two triangles glued together along the edges
corresponding to the amalgamated actions:
\bea
\source&=&\sourcepic\\
\target&=&\targetpic
\eea
The numbering of the faces refer to their order in the tensor
product. An other suggestive notation would be to denote the
domain of $V$ by $\Hil\oR\Hil$ and its range by $\Hil\oL\Hil$.
We can now define the operator
$U\colon\Hil\oR\Hil\to\Hil\oL\Hil$ as the
restriction of $V$ to its domain and range. The natural
representation of this operator is then the tetrahedron
\beq
U\ =\ \tetra
\eeq
or, better to say, this projection of the tetrahedron. Namely,
the "equator" $\{(01),(12),(23),(30)\}$ is distinguished by
dividing the surface into a "Northern hemisphere"
$\{(012),(023)\}$ and a "Southern hemisphere"
$\{(013),(123)\}$ corresponding
to the range and domain of $U$, respectively.
 
Both the range and domain of $U$ are quadrimodules, i.e. $A^L$
acts on them via 3 left actions
$\alpha_{01},\alpha_{12},\alpha_{23}$ and 1 right action
$\alpha_{03}$, and these 4 perimeter actions commute with each
other. For example $\alpha_{12}$ acts on $\Hil\oR\Hil$ as
$\id\o\alpha_{01}$ and on $\Hil\oL\Hil$ as $\alpha_{12}\o\id$.
Now $U$ can be shown to intertwine these four actions,
$\alpha_{ij}(x^L)U=U\alpha_{ij}(x^L)$, $x^L\in A^L$,
$(ij)=(01),(12),(23),(03)$. The intertwiner relations are
consequences of the following identities for $V$:
\beabc  \label{inter}
V(x^L\o\one)&\stackrel{(01)}{=}&(x^L\o\one)V\qquad x^L\in A^L\\
V(\one\o x^L)&\stackrel{(12)}{=}&(\varphi^R\o\one)V\qquad
\varphi^R=x^L\la\du1,\ x^L\in A^L\\
V(\one\o\varphi^R)&\stackrel{(23)}{=}&(\one\o\varphi^R)V\qquad
\varphi^R\in\duA^R\\
V(x^R\o\one)&\stackrel{(03)}{=}&(\one\o x^R)V\qquad x^R\in A^R
\eeabc
(Here $\stackrel{(ij)}{=}$ refers to the edge $(ij)$ of the
tetrahedron $(0123)$ and not to an equation number as before.)
The intertwiner relations for $U$ are precisely the four equations
in Definition 5.6.i of \cite{EV}. Thus, in order to see that
our $U$ is a pseudo-multiplicative unitary,
we are left with showing that $U$ satisfies the pentagon equation
of Figure 1.
\begin{figure}
\begin{picture}(400,460)(0,-20)
\put(0,400){$\t{0}{1}{4}\amalgo{14}\t{1}{2}{4}\amalgo{24}\t{2}{3}{4}$}
\put(50,360){$1\t{0}{1}{4}\amalgo{14} U(^{14}_{23})$}
\put(0,320){$\t{0}{1}{4}\amalgo{14}[\t{1}{2}{3}\amalgo{13}\t{1}{3}{4}]$}
\put(85,280){$1\t{0}{1}{4}\amalgo{14}\Sigma(\two,\three)$}
\put(105,240){$\t{0}{1}{4}\amalgo{14}\t{1}{3}{4}\amalgo{13}\t{1}{2}{3}$}
\put(95,200){$\Sigma(\one\two,\three)$}
\put(0,160){$\t{1}{2}{3}\amalgo{13}[\t{0}{1}{4}\amalgo{14}\t{1}{3}{4}]$}
\put(50,120){$1\t{1}{2}{3}\amalgo{13} U(^{04}_{13})$}
\put(0,80){$\t{1}{2}{3}\amalgo{13}\t{0}{1}{3}\amalgo{03}\t{0}{3}{4}$}
\put(85,40){$\Sigma(\one,\two)\amalgo{03}1\t{0}{3}{4}$}
\put(105,0){$[\t{0}{1}{3}\amalgo{13}\t{1}{2}{3}]\amalgo{03}\t{0}{3}{4}$}
 
\put(45,385){\vector(0,-1){50}}
\put(45,305){\vector(1,-1){50}}
\put(95,225){\vector(-1,-1){50}}
\put(45,145){\vector(0,-1){50}}
\put(45,65){\vector(1,-1){50}}
 
\put(93,425){$U(^{04}_{12})\amalgo{24} 1\t{2}{3}{4}$}
\put(170,400){$\t{0}{1}{2}\amalgo{02}\t{0}{2}{4}\amalgo{24}\t{2}{3}{4}$}
\put(261,425){$1\t{0}{1}{2}\amalgo{02} U(^{04}_{23})$}
\put(340,400){$\t{0}{1}{2}\amalgo{02}\t{0}{2}{3}\amalgo{03}\t{0}{3}{4}$}
 
\put(90,405){\vector(1,0){70}}
\put(260,405){\vector(1,0){70}}
 
\put(280,200){$U(^{03}_{12})\amalgo{03} 1\t{0}{3}{4}$}
\put(205,10){\vector(1,3){125}}
\end{picture}
\caption{The `pentagon' equation for the pseudo-multiplicative
unitary $U$.}
\end{figure}
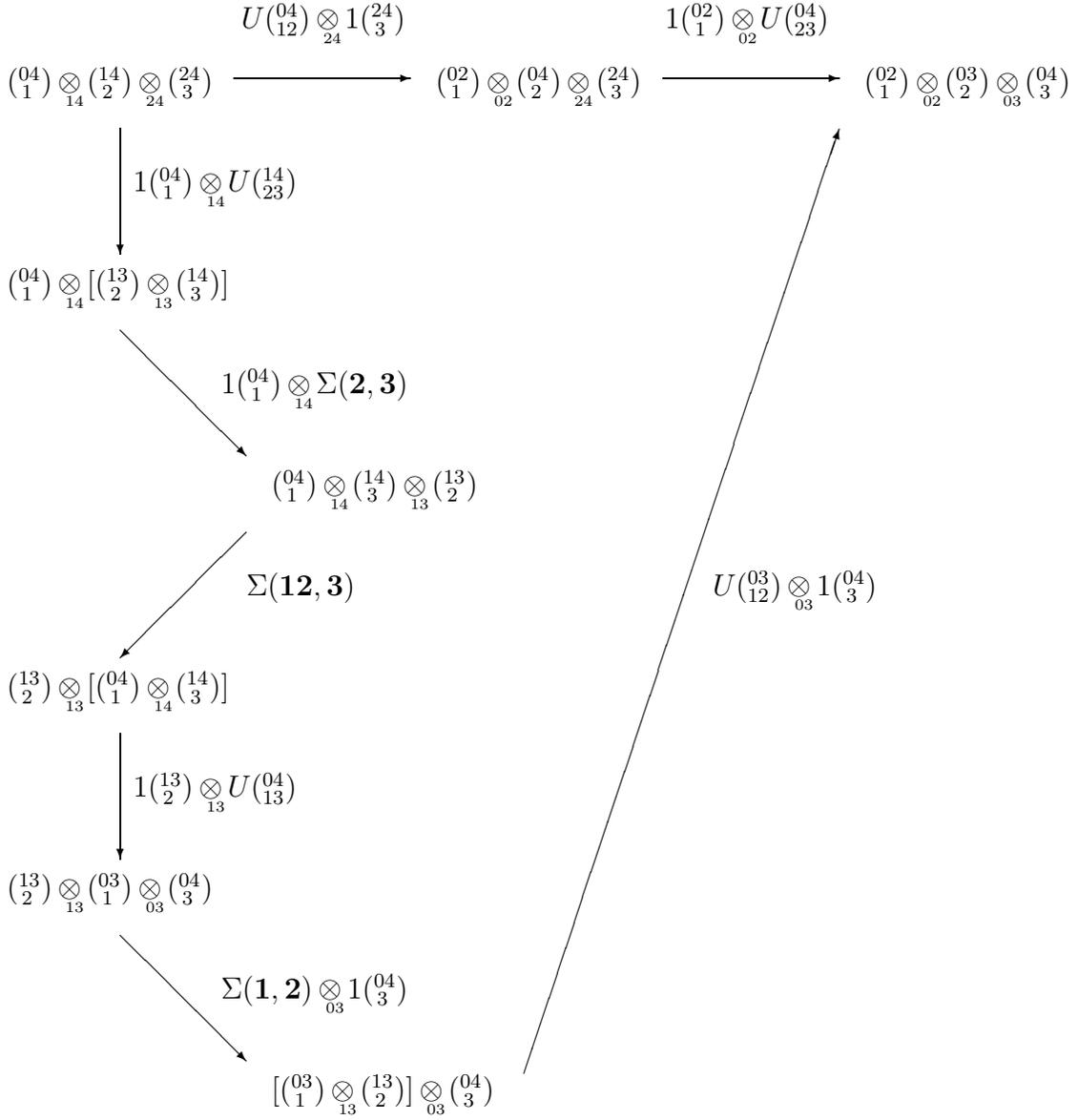
Before doing that we remark that
on the remaining two edges of the tetrahedron we have the
amalgamation relations (\ref{heart2}) and (\ref{heart1}),
\bea \label{ams}
\xi\amalgo{\varepsilon^o} \alpha_{02}(x^L)\eta &=&
\alpha_{12}(g_L^{\half} x^L g_L^{-\half})\xi\amalgo{\varepsilon^o}\eta
\\
\label{amt}
\xi\amalgo{\varepsilon} \alpha_{01}(x^L)\eta &=&
\alpha_{02}(g_L^{-\half}x^L g_L^{\half})\xi\amalgo{\varepsilon} \eta\ ,
\eea
which have their origin in the $V$-identities
\bea
V(\varphi^R\o\one)&=&V(\one\o x^R)\qquad \varphi^R=x^R\la\du1\\
(\varphi^L\o\one)V&=&(\one\o x^L)V\qquad \varphi^L=\du1\ra x^L
\eea

As for the pentagon equation is concerned we need a more concise
notation for multiple relative tensor products. Therefore we use
the symbol $\t{i}{j}{k}$ to denote a copy of $\Hil$ associated
to the triangle $(ijk)$. The symbol $\amalgo{ij}$ will stand for
the relative tensor product of the two triangle modules that
contain the edge $(ij)$. Whether it is a tensor product with
respect to $\eps$ or $\eps^o$ can be unambiguously recovered
from the order of the modules in the tensor product. This is
because  each internal edge $(ij)$
(of a planar 2-complex) has opposite relative orientation to its
two neighbour faces. For example
$$
(^{02}_{\,1})\amalgo{02}(^{04}_{\,2})\amalgo{24}(^{24}_{\,3})\ =\
\Hil_{\alpha_{02}}\amalgo{\eps}\ _{\alpha_{01}}\Hil_{\alpha_{12}}
\amalgo{\eps^o}\ _{\alpha_{02}}\Hil
$$
Sometimes it is unavoidable to use brackets because $\amalgo{ij}$
refers to two triangles that are not consecutive ones in the
tensor product. For example in
$$
[(^{03}_{\,1})\amalgo{13}(^{13}_{\,2})]\amalgo{03}(^{04}_{\,3})\ =
\ [\Hil_{\alpha_{12}}\amalgo{\eps^o}\ _{\alpha_{02}}\Hil]
_{\alpha_{03}}\amalgo{\eps}\ _{\alpha_{01}}\Hil
$$
These brackets therefore have nothing to do with associativity of
the tensor product. They reflect rather the poor capability of
our one dimensional writing to express two dimensional facts.
 
Now we are ready to formulate the pentagon equation.
In our notation the equation of Definition 5.6.ii of \cite{EV}
takes the form as Figure 1.
The boldface numbers in the argument of the flip map
refer to factors of the tensor product that forms the domain of
$\Sigma$. E.g. $\Sigma(\one\two,\three)$ maps $\xi\o\eta\o\zeta$
to $\zeta\o\xi\o\eta$. Up to the flip maps, which serve only for
permuting the tensor product factors in linear writing, the above
commutative diagram is a pentagon rather than an octogon.
The reader may find it amusing to draw the eight pentagonal
figures corresponding to the eight vertices of Figure 1, each of
them having vertices numbered (01234), have two diagonals
one for each $\amalgo{ij}$ symbol, and have triangular faces
numbered according to their order in the tensor product.
Let stand here one edge of Figure 1 for example:
 
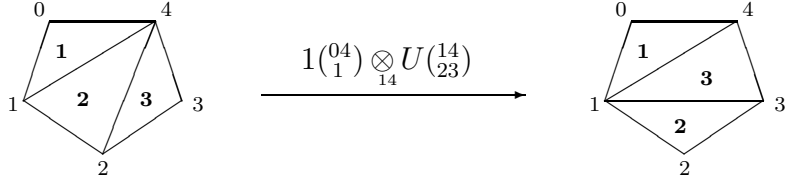
\begin{picture}(290,100)(-20,0)
\put(0,40){\penta}
\put(100,40){\vector(1,0){100}}
\put(115,50){$1\t{0}{1}{4}\amalgo{14}U(^{14}_{23})$}
\put(220,40){\pentta}
\end{picture}
 
After acquainting the equation we have to show that it is a
consequence of the $V$-pentagon (\ref{penta}). At first we
identify the eight corners in Figure 1 with subspaces of
$\Hil\o\Hil\o\Hil$. With the notation $E=VV^*$, $\hat E=V^*V$
we can write
\beanon
\t{0}{1}{4}\amalgo{14}\t{1}{2}{4}\amalgo{24}\t{2}{3}{4}&=&
\hat E_{12}\hat E_{23}(\Hil\o\Hil\o\Hil)\\
\t{0}{1}{4}\amalgo{14}[\t{1}{2}{3}\amalgo{13}\t{1}{3}{4}]&=&
\hat E_{13} E_{23}(\Hil\o\Hil\o\Hil)\\
\t{0}{1}{4}\amalgo{14}\t{1}{3}{4}\amalgo{13}\t{1}{2}{3}&=&
\hat E_{12} E_{32}(\Hil\o\Hil\o\Hil)\\
\t{1}{2}{3}\amalgo{13}[\t{0}{1}{4}\amalgo{14}\t{1}{3}{4}]&=&
 E_{13}\hat E_{23}(\Hil\o\Hil\o\Hil)\\
\t{1}{2}{3}\amalgo{13}\t{0}{1}{3}\amalgo{03}\t{0}{3}{4}&=&
\hat E_{21} E_{23}(\Hil\o\Hil\o\Hil)\\
{[}\t{0}{1}{3}\amalgo{13}\t{1}{2}{3}]\amalgo{03}\t{0}{3}{4}&=&
 E_{13}\hat E_{12}(\Hil\o\Hil\o\Hil)\\
\t{0}{1}{2}\amalgo{02}\t{0}{2}{4}\amalgo{24}\t{2}{3}{4}&=&
 E_{12}\hat E_{23}(\Hil\o\Hil\o\Hil)\\
\t{0}{1}{2}\amalgo{02}\t{0}{2}{3}\amalgo{03}\t{0}{3}{4}&=&
 E_{12} E_{23}(\Hil\o\Hil\o\Hil)
\eeanon
Inserting $V=VV^*V=EV\hat E$ into the $V$-pentagon
(\ref{penta}) we obtain
\[
E_{23}V_{23}(\hat E_{23}E_{12})V_{12}\hat E_{12}=
E_{12}V_{12}(\hat E_{12}E_{13})V_{13}(\hat E_{13}E_{23})V_{23}\hat
E_{23}\ .
\]
Multiplying with projections from the left and right
and inserting appropriate flip maps
\bea
(E_{12}E_{23})V_{23}(\hat E_{23}E_{12})V_{12}(\hat E_{12}\hat
E_{23})&=&(E_{12}E_{23})V_{12}(\hat E_{12}E_{13})\Sigma_{1,2}
(\hat E_{21}E_{23})\nn
&&V_{23}(\hat E_{23}E_{13})\Sigma_{12,3}(\hat E_{12}E_{32})
\Sigma_{2,3}(\hat E_{13}E_{23})\nn
&&V_{23}(\hat E_{12}\hat E_{23})\label{pentaindisguise}
\eea
The eight different projections in the parentheses correspond
precisely to the eight corners of the diagram in Fig.1. The $V$
and $\Sigma$ operators, together with their neighbour projections,
in turn produce precisely the eight maps of the diagram. In order
to see this one should check correspondences like
\beanon
U(^{04}_{12})\amalgo{24}1\t{2}{3}{4}&\equiv&\hat E_{23}V_{12}\hat
E_{23}=(E_{12}\hat E_{23})V_{12}(\hat E_{12}\hat E_{23})\\
1\t{0}{1}{2}\amalgo{02} U(^{04}_{23})&\equiv&E_{12}V_{23}E_{12}=
(E_{12}E_{23})V_{23}(E_{12}\hat E_{23})\\
1\t{0}{1}{4}\amalgo{14}\Sigma(\t{1}{2}{3},\t{1}{3}{4})&\equiv&
(\hat E_{12}E_{32})\Sigma_{23}(E_{23}\hat E_{13})
\eeanon
and five other ones. This finishes the proof of that every unital
MPI $V$ determines a pseudo-multiplicative unitary $U$ by
restriction to range and domain. As a byproduct we obtained a
geometric interpretation of the equations in terms of trimodules,
or 2-simplex modules, $\Hil$ over $A^L$ in which $U$ plays the
role of Ocneanu's 3-cocycle.
 
\medskip
Now we turn to the opposite construction when we are given a
pseudo-multiplicative unitary $U$
and want to construct a multiplicative partial isometry
$V$ that reproduces $U$ by restriction. This task will be a simple
one mainly
because we can prove only that the resulting $V$ is an MPI and we
leave it open whether $V$ is unital.
 
Let $N$ be a finite dimensional $C^*$-algebra with a faithful
positive linear functional $\nu\colon N\to\C$ of index $1$. Let
$\beta,\alpha,\hat\beta$ be actions of $N$, $N^o$, and $N$,
respectively on a finite dimensional Hilbert space $\Hil$ that
commute with each other. Finally,
let $U\colon\Hil_{\hat\beta}\amalgo{\nu^o}\,_{\alpha}\Hil\to
\Hil_{\alpha}\amalgo{\nu}\,_{\beta}\Hil$ be a
pseudo-multiplicative isometry.
 
Since the relative tensor products can be identified as subspaces
in $\Hil\o\Hil$ via the projections (\ref{e}), we can immediately
define a partial isometry
\beq
V:=EU\hat E\,,\quad\mbox{where}\quad E=(\alpha\o\beta)(e_\nu),\
\hat E=(\hat\beta\o\alpha)(e_\nu)\,.
\eeq
Then $VV^*=E$ and $V^*V=\hat E$. Defining the algebras
$A^L:=\beta(N)$, $A^R:=\alpha(N)$, and $\duA^R:=\hat\beta(N)$
the four intertwiner relations for $U$ become the intertwiner
relations (\ref{inter}-d). These in turn are equivalent to the
equations (\ref{hexa2}), (\ref{hexa1}), (\ref{hexa3}), and
(\ref{fifth}), respectively. The pentagon equation (\ref{penta})
can now be obtained by arguing backwards with equation
(\ref{pentaindisguise}). This proves that $V$ is a multiplicative
isometry.

\end{document}